\documentclass[10pt]{amsproc}
\usepackage{graphicx}
\usepackage{amscd}
\usepackage{amsmath}
\usepackage{amsfonts}
\usepackage{amssymb}
\newtheorem{theorem}{Theorem}[section]
\newtheorem{corollary}[theorem]{Corollary}
\newtheorem{lemma}[theorem]{Lemma}
\newtheorem{proposition}[theorem]{Proposition}
\theoremstyle{definition}
\newtheorem{definition}[theorem]{Definition}

\theoremstyle{remark}
\newtheorem{remark}[theorem]{Remark}

\numberwithin{equation}{section}

\begin{document}
\title[Invariant Subspace Theorem for Multivalued Operators]{Lomonosov's Invariant Subspace Theorem for Multivalued Linear Operators}
\author{Peter Saveliev}
\address{Allegheny College, Meadville, PA 16335}
\email{saveliev@member.ams.org}
\thanks{}
\date{May 11, 2000}

\begin{abstract}
The famous Lomonosov's invariant subspace theorem states that \textit{if a
continuous linear operator }$T$\textit{\ on an infinite-dimensional normed
space }$E$\textit{\ ``commutes'' with a compact operator }$K\neq0,$ i.e.,
$TK=KT,$\textit{\ then }$T$\textit{\ has a non-trivial closed invariant
subspace}. We generalize this theorem for multivalued linear operators. We
also provide some applications to single-valued linear operators.
\end{abstract}\maketitle

\section{Introduction.}

The Invariant Subspace Problem asks whether every linear operator
$h:Y\rightarrow Y$ on an infinite dimensional topological vector space $Y$ has
a nontrivial closed \textit{invariant subspace}, i.e., a linear subspace $M$
of $Y$ such that $M\neq\{0\},$ $M\neq Y$ and $h(M)\subset M$ (for a survey and
references see \cite{RR,AAB}). In general the answer is negative and the issue
is to investigate the class of operators satisfying this property. It is known
that every compact operator $k$ belongs to this class and so does every
operator $h$ commuting with $k:$
\[
hk=kh.
\]
The famous Lomonosov's invariant subspace theorem \cite{Lom} states the following.

\begin{theorem}
\label{LomOrig}Suppose $Y$ is an infinite dimensional normed space, $h$ and
$k$ are continuous linear operators, $k$ is compact, nonzero and commutes with
$h$. Then $h$ has a nontrivial closed invariant subspace.
\end{theorem}

In this paper we generalize this result for multivalued linear operators. The
theory of multivalued linear maps (linear relations) is well developed, see
Cross \cite{Cross} A multivalued map $h:X\rightarrow Y$ between vector spaces
is called a \textit{linear relation} if
\[
h(ax)=ah(x),\quad h(x+y)=h(x)+h(y),
\]
for all $x,y\in X$ and all $a\neq0.$ We say that $h,k:Y\rightarrow Y$
\textit{commute} if
\[
hk\subset kh.
\]
For a linear relation $h:Y\rightarrow Y,$ a subspace $M$ of $Y$ is called
$h$\textit{-weakly-invariant} if for all $x\in M$%
\[
h(x)\cap M\neq\emptyset.
\]
The idea of the proof of our main theorem below (Corollary \ref{InvSubCor})
can be traced back to the original Lomonosov's proof.

\begin{theorem}
Suppose $Y$ is an infinite dimensional normed topological vector space, $h$
and $k$ are continuous linear relations with nonempty finite dimensional
values, $k$ is compact and commutes with $h$, $k^{-1}(0)\neq Y.$ Then there is
a nontrivial closed $h$-weakly-invariant subspace.
\end{theorem}

Invariant subspace theorems for linear relations provide tools for studying
the Invariant Subspace Problem for single valued linear operators. We consider
those in the last section.

All topological spaces are assumed to be Hausdorff, all maps are multivalued
with \textit{nonempty} values unless indicated otherwise, by normed (locally
convex) spaces we understand infinite dimensional normed (locally convex)
topological vector spaces over $\mathbf{C}$ or $\mathbf{R}$.\ \ 

\section{Preliminaries.}

Let $X$ be a topological space. The a \textit{partition of unity} is a
collection of continuous functions \ $\gamma=\{d_{\alpha}:\alpha\in A\}$
satisfying
\[%
{\displaystyle\sum\limits_{\alpha\in A}}
d_{\alpha}(x)=1,x\in X.
\]
The partition $\gamma$ is called \textit{locally finite} if the cover
$\gamma^{\prime}=\{d_{\alpha}^{-1}((0,1]):\alpha\in A\}$ of $X$ is locally
finite, and $\gamma$ is called \textit{subordinate to an open cover} $\omega$
of $X$ if $\gamma^{\prime}$ refines $\omega.$

Let $F:X\rightarrow Y$ be a multifunction (a set-valued map $F:X\rightarrow
2^{Y}$), where $X,Y$ are topological spaces. We call $F$
\textit{lower-semicontinuous }(l.s.c.) if $F^{-1}(B)=\{x\in X:F(x)\cap
B\neq\emptyset\}$ is open for any open $B$. We call $F$
\textit{upper-semicontinuous }(u.s.c.) if $F^{-1}(B)$ is closed for any closed
$B$. Equivalently, for any $x\in X$ and an neighborhood $V$\ of $F(x),$ there
is a neighborhood $U$ of $x$ such that
\[
y\in U\Longrightarrow F(y)\subset V.
\]
If $F$ is both u.s.c. and l.s.c. we call it \textit{continuous} (in
\cite{Cross} ``continuous'' means l.s.c.). When $X$ and $Y$ are uniform spaces
\cite[Chapter 8]{En} we say that $F$ is \textit{uniformly}
\textit{upper-semicontinuous} (u.u.s.c.) if for any entourage $V$ in $Y$ there
is an entourage $U$ in $X,$ such that for all $x,y\in X,$%
\[
x\in y+U\Longrightarrow F(x)\subset F(y)+V.
\]
We say that $x\in X$ is a \textit{fixed point} of $F:X\rightarrow X$ if $x\in
F(x).$ We say that a single-valued function $g:X\rightarrow Y$ is a
\textit{selection} of $F:X\rightarrow Y$ if $g(x)\in F(x)$ for all $x\in X$.

A topological space is called \textit{acyclic} if its reduced \v{C}ech
homology groups over the rationals are trivial. The set of all acyclic spaces
includes convex subsets of locally convex spaces. A multivalued map
$F:X\rightarrow Y$ is called \textit{admissible} (in the sense of Gorniewicz
\cite{Gorn}) if it is closed valued u.s.c. and there exist a topological space
$Z$ and two single-valued continuous maps $p:Z\rightarrow X,\ q:Z\rightarrow
Y$ such that for any $x\in X,$ (i) $p^{-1}(x)$ is acyclic, and (ii)
$q(p^{-1}(x))\subset F(x).$

\begin{proposition}
\label{AdmComp}\cite[Theorem IV.40.5, p. 200]{Gorn} The composition of two
admissible maps is admissible.
\end{proposition}

In the proof of his theorem Lomonosov used the Schauder fixed point theorem.
Here we need its analogue for multivalued maps. We use a result that follows
from the Lefschetz fixed point theory for admissible maps given in \cite{Gorn}
(for further developments see \cite{Sav1}), although an appropriate version of
the Kakutani fixed point theorem for compositions of convex valued maps would suffice.

\begin{theorem}
\label{GornFP}\cite[Theorem IV.41.12, p. 207]{Gorn} If $X$ is an acycic ANR,
then $X$ has the fixed point property within the class of compact (i.e., with
$\overline{f(X)}$ compact) admissible maps.
\end{theorem}

\section{Approximation by Convex Combinations.}

In this section we consider the approximations of multivalued maps by convex
combinations in the spirit of Simoni\v{c} \cite{Sim1,Sim2}.

Let $X$ be a topological space, $Y$ be a vector space and $H$ a collection of
multivalued maps $h:X\rightarrow Y.$ Given a partition of unity$\ \gamma
=\{d_{\alpha}:\alpha\in A\}$ in $X,$ for any collection $\{h_{\alpha}%
:\alpha\in A\}\subset H,$ we can define a new map $h:X\rightarrow Y$ by
\begin{equation}
h(x)=%
{\displaystyle\sum\limits_{\alpha\in A}}
d_{\alpha}(x)h_{\alpha}(x). \label{conv}%
\end{equation}
The set of all such $h$ we denote by $Conv_{\gamma}(H).$

Let $\Delta_{n}$ denote the standard $n$-simplex. For any collection
$h_{0},...,h_{n}$ and any $(t_{0},...,t_{n})\in\Delta_{n},$ we define a new
map $h:X\rightarrow Y$ by
\[
h(x)=%
{\displaystyle\sum\limits_{i=0}^{n}}
t_{i}h_{i}(x).
\]
The set of all such $h$ we denote by $Conv(H).$

\begin{remark}
$H\subset Conv(H)\subset Conv_{\gamma}(H).$
\end{remark}

The following theorem generalizes Lemma 3.1 of Simoni\v{c} \cite{Sim1}.

\begin{theorem}
[Approximation]\label{Appr}Let $X$ be a paracompact uniform space and $Y$ be a
locally convex space. Suppose $s:X\rightarrow Y$ is u.u.s.c. map with convex
values and $H$ is a collection of u.s.c. maps $h:X\rightarrow Y.$ Suppose $V$
is a convex neighborhood of $0$ in $Y$ and for any $x\in X$ there is $h_{x}\in
H$ such that
\[
h_{x}(x)\subset s(x)+V.
\]
Then there exist a locally finite partition of unity $\gamma$ on $X$ and a map
$f\in Conv_{\gamma}(H)$ with values such that for all $x\in X,$%
\[
f(x)\subset s(x)+3V.
\]
\end{theorem}

\begin{proof}
[Proof]From the upper semicontinuity it follows that for each $x\in X$ there
is a $W_{x}$ such that for any $y\in W_{x},$%
\[
h_{x}(y)\subset h_{x}(x)+V,\text{ and }s(x)\subset s(y)+V.
\]
Now applying these inclusions and the assumption of the theorem, we obtain the
following: for any $y\in W_{x},$
\[
h_{x}(y)\subset h_{x}(x)+V\subset s(x)+2V\subset s(y)+3V.
\]

Consider the open cover of $X$ given by $\omega=\{W_{x}:x\in X\}.$ From
Michael's Lemma \cite[Theorem 5.1.9, p. 301]{En} it follows that there exists
a locally finite partition of unity $\gamma=\{d_{\alpha}:\alpha\in A\}$
subordinate to $\omega.$ Then we have a locally finite open cover of $X $%
\[
\omega^{\prime}=\{U_{\alpha}=d_{\alpha}^{-1}((0,1]):\alpha\in A\}
\]
that refines $\omega,$ i.e., for each $\alpha\in A$ there is $x(\alpha)\in X$
such that $U_{\alpha}\subset W_{x(\alpha)}.$

Fix $y\in X.$ Suppose $\alpha\in A$ and $d_{\alpha}(y)>0.$ Then $y\in
U_{\alpha}\subset W_{x(\alpha)}.$ Hence
\[
h_{x(\alpha)}(y)\subset s(y)+3V.
\]
As $s(y)$ and $V$ are convex, so is the set $s(y)+3V$. Therefore a convex
combination of $h_{x(\alpha)}(y),\alpha\in A$ is a well defined subset of
$s(y)+3V$. Then the map $f:X\rightarrow Y$ given by
\[
f(y)=%
{\displaystyle\sum\limits_{\alpha\in A}}
d_{\alpha}(y)h_{\alpha}(y)
\]
is well defined and belongs to $Conv_{\gamma}(H).$
\end{proof}

Simoni\v{c} calls the functions $d_{\alpha}$ \textit{Lomonosov functions} as
the idea of this construction goes back to Lomonosov's proof in \cite{Lom}.
Observe also that this theorem implies the following well-known fact:
\textit{any u.s.c. map with convex images from a compact metric space to a
locally convex space can be approximated by continuous single-valued maps}.

\section{Fixed Points of Convex Combinations.}

In this section we obtain a preliminary fixed point result.

\begin{proposition}
\label{Fix}Let X be a topological space, $Y$ a vector space, $H$ be a
collection of maps $h:X\rightarrow Y,$ $r:Y\rightarrow X$ a map$,$ $\gamma$ a
locally finite partition of unity on $X,$ $f\in Conv_{\gamma}(H),$ and suppose
$fr:Y\rightarrow Y$ has a fixed point $y_{0}.$ Then there exists $g\in
Conv(H)$ such that $y_{0}$ is a fixed point of $gr.$
\end{proposition}

\begin{proof}
We know that $y_{0}\in fr(y_{0}).$ Let $Z=Graph(r)\subset Y\times X,$
$p:Z\rightarrow y,q:Z\rightarrow X$ be the projections. Then there \ is
$z_{0}\in Z$ such that $y_{0}=p(z_{0})\in fq(z_{0}).$ Assume that
$\gamma=\{d_{\alpha}:\alpha\in A\}$ and suppose
\[
f(x)=%
{\displaystyle\sum\limits_{\alpha\in A}}
d_{\alpha}(x)h_{\alpha}(x),
\]
where $h_{\alpha}\in H.$ Let $x_{0}=q(z_{0}).$ Suppose $\{\alpha\in
A:d_{\alpha}(x_{0})>0\}=\{\alpha_{0},...,\alpha_{n}\},n\geq0.$ For
$i=0,...,n,$ let
\[
t_{i}=d_{\alpha_{i}}(x_{0}),h_{i}=h_{\alpha_{i}}.
\]
Then $y_{0}\in fq(z_{0})=f(x_{0})=\sum_{i=0}^{n}t_{i}h_{i}(x_{0}).$ We have an
element of $Conv(H):$%
\[
g(x)=\sum_{i=0}^{n}t_{i}h_{i}(x),x\in X.
\]
Consider
\[%
\begin{array}
[c]{ll}%
gq(z_{0}) & =\sum_{i=0}^{n}t_{i}h_{i}(q(z_{0}))\\
& =\sum_{i=0}^{n}d_{\alpha_{i}}(x_{0})h_{i}(q(z_{0}))\\
& =\sum_{\alpha\in A}d_{\alpha_{i}}(q(z_{0}))h_{i}(q(z_{0}))\\
& =fq(z_{0})\\
& \ni p(z_{0})=y_{0}.
\end{array}
\]
Thus $y_{0}$ is a fixed point of $gr.$
\end{proof}

If $Y$ is a topological vector space, we denote by $F_{c}(X,Y)$ the set of all
u.s.c. maps $F:X\rightarrow Y$ with compact convex values. Then all elements
of $F_{c}(X,Y)$ are admissible provided $Y$ is locally convex.

\begin{lemma}
\label{Fconv}$Conv_{\gamma}(F_{c}(X,Y))\subset F_{c}(X,Y).$
\end{lemma}

\begin{theorem}
[Fixed Points]\label{FixNEmp}Let $Y$ be a locally convex space, $A\subset Y$
closed convex, $U$ a closed neighborhood of $A,$ $X$ a paracompact uniform
space. Let $r\in F_{c}(U,X)$, $\overline{r(U)}$ compact, $H\subset
F_{c}(X,Y),$ and
\[
\overline{H(x)}\cap A\neq\emptyset,\text{ for all }x\in X.
\]
Then for any neighborhood $W$ of $A$, $Conv(H)r$ has a fixed point in $W.$
\end{theorem}

\begin{proof}
Assume that there is a convex neighborhood $V$ of $0$ such that $U=A+3V\subset
W$ and for each $x\in X,$ there is $h_{x}\in H$ such that $h_{x}%
(x)\cap(A+1/2V)\neq\emptyset.$ Now, for each $h\in H,$ let
\[
h^{\prime}(x)=h(x)\cap\overline{(A+1/2V)}%
\]
and let $H^{\prime}=\{h^{\prime}:h\in H\}.$ Then the set $h^{\prime}(x)$ is
nonempty by assumption, convex since $Y$ is locally convex, compact as the
intersection of a compact set and a closed set. Also $h^{\prime}$ is u.s.c. by
\cite[1.7.17(c)]{En}. Thus $h^{\prime}\in F_{c}(X,Y).$ Now we apply Theorem
\ref{Appr} with $s(x)=A$ for all $x$ (by definition $h_{x}^{\prime}(x)\subset
s(x)+V).$ Therefore there exists $f\in Conv_{\gamma}(H^{\prime})$ such that
for all $x,$ $f(x)\subset s(x)+3V=U.$ We know that $r,f$ are admissible.
Therefore by Proposition \ref{AdmComp}, so is $\varphi=fr:U\rightarrow U$ .
Now $r$ is compact and $U$ is an ANR as a closed neighborhood in a locally
convex space. Hence by Theorem \ref{GornFP}, $\varphi$ has a fixed point
$y_{0}\in U\subset W.$ Therefore by Proposition \ref{Fix}, there exists $g\in
Conv(H)$ such that $y_{0}$ is a fixed point of $gr.$
\end{proof}

\section{Properties of Linear Relations.}

Throughout the rest of the paper we assume that $X$ and $Y$ are normed spaces.

\begin{definition}
\cite[II.1.3, p. 25]{Cross} A multivalued map $h:X\rightarrow Y$ is called a
\textit{linear relation} if it preserves nonzero linear combinations, i.e.,
for all $x,y\in X,$ all $a,b\in\mathbf{R}\backslash\{0\},$ we have
\[
h(ax+by)=ah(x)+bh(y).
\]
Then $h(0)$ is a linear subspace. The set of all linear relations will be
denoted by $LR(X,Y),$ and $LR(X,X)=LR(X).$
\end{definition}

\begin{lemma}
\label{T(0)}\cite[Proposition I.2.8, p. 7]{Cross} If $T\in LR(X,Y)$, $x\in X,
$ then
\[
T(x)=y+T(0),\text{ for any }y\in T(x).
\]
\end{lemma}

We will concentrate on the following classes of linear relations:
\begin{align*}
LR_{0}(X,Y)  &  =\{h\in LR(X,Y):h\text{ is continuous, }\dim h(0)<\infty\},\\
LR_{0}(Y)  &  =LR_{0}(Y,Y).
\end{align*}
\emph{\ }Of course, all bounded linear operators belong to $LR_{0}(X,Y).$

For a linear relation $T\in LR_{0}(X,Y),$ let $Q_{T}$ denote the natural
quotient map with domain $X$ and null space $T(0)$ \cite[p. 25]{Cross}.

\begin{lemma}
If $S\in LR_{0}(X,Y),$ $T\in LR_{0}(Y,Z),$ then $TS\in LR_{0}(X,Z).$
\end{lemma}

\begin{proof}
$Q_{T}T:S(0)\rightarrow Z/T(0)$ is a linear operator, so $\dim Q_{T}%
TS(0)<\infty$. Now $\dim T(0)<\infty$ implies $\dim TS(0)<\infty.$
\end{proof}

\begin{definition}
\cite{Cross} A linear relation $T\in LR_{0}(X,Y)$ is called \textit{bounded}
(\textit{compact}) if the single valued operator $Q_{T}T$ is bounded
(compact), i.e., it maps a bounded set into a bounded (compact) set.
\end{definition}

By Proposition II.3.2(a) in \cite[p. 33]{Cross}, every element of
$LR_{0}(X,Y)$ is bounded.

\begin{lemma}
\label{Cmpct}Suppose $Y$ is normed. Then $T\in LR_{0}(X,Y)$ is compact if and
only if for any bounded $B\subset X,$ there is a compact set $C\subset Y$ such
that $T(B)\subset C+T(0).$ Moreover $C$ can be chosen such that $T(x)\cap
C\neq\emptyset$ for all $x\in B$.
\end{lemma}

\begin{proof}
The ``if'' part is obvious. Next, if $Q_{T}T:X\rightarrow Z=Y/T(0)$ is a
compact linear operator, then for any bounded $B\subset X,$ there is a compact
$D\subset Z$ such that $Q_{T}T(B)\subset D.$ Now since $T(0)$ is a finite
dimensional subspace of a normed space, it is topologically complemented,
i.e., $Y$ is homeomorphic to $Z\oplus T(0).$ Then $C=D\oplus\{0\}$ is compact
in $Y$ and $T(B)\subset C+T(0)$.
\end{proof}

\begin{theorem}
\label{CompactComposition}Suppose $h\in LR_{0}(X,Y)$ and $k\in LR_{0}(Y,Z)$ is
compact. Then $hk\in LR_{0}(X,Z)$ is compact.
\end{theorem}

\begin{proof}
Let $A$ be a bounded subset of $X.$ Then by Lemma \ref{Cmpct}, $k(A)\subset
C+k(0),$ where $C$ is compact. Therefore $hk(A)\subset h(C)+hk(0).$ It is easy
to show that there is a bounded set $D$ in $Z$ such that $h(C)=D+h(0).$ Let
$h^{\prime}(x)=h(x)\cap\overline{D}.$ Then $h^{\prime}(x)$ is compact as
$h(x)$ is finite dimensional. In particular, $h^{\prime}(x)$ is closed, so by
\cite[1.7.17(c)]{En}, $h^{\prime}=h\cap\overline{D}$ is u.s.c.. Therefore
$h(C)\cap\overline{D}=h^{\prime}(C)$ is compact in $Z$ by \cite[Proposition
II.14.9, p. 69]{Gorn}. But $D\subset h^{\prime}(C),$ hence $D$ is precompact.
Thus
\[
hk(A)\subset h(C)+hk(0)=D+h(0)+hk(0)=D+hk(0),
\]
where $D$ is precompact, so $hk$ is compact by Lemma \ref{Cmpct}.
\end{proof}

Therefore compact relations constitute a left ideal in $LR_{0}(Y).$

We call $G\subset LR(Y)$ a \textit{semialgebra} if it is closed under nonzero
linear combinations and compositions. The above lemma implies that $LR_{0}(Y)$
is a semialgebra. We define the \textit{commutant} of $h$ as
\[
Comm(h)=\{f\in LR_{0}(Y):f\text{ commutes with }h\}.
\]

\begin{lemma}
\label{CommLin} If $h\in LR(Y)$ then $Comm(h)$ is a semialgebra.
\end{lemma}

\begin{proof}
First, $Comm(h)$ is closed under linear combinations. Indeed, for
$f,f^{\prime}\in Comm(h),$ $a,b\in\mathbf{R\backslash\{0\}},$ we have
\[
h(af+bf^{\prime})=ahf+bhf^{\prime}\subset afh+bf^{\prime}h=(af+bf^{\prime})h.
\]
Second, if $f,g$ both commute with $h,$ then $u=gf$ commutes with $h.$ Indeed
\[
hu=hgf\subset ghf\subset gfh=uh.
\]
\end{proof}

\section{Invariant Subspaces of Linear Relations.}

\begin{lemma}
\label{ComLin}Suppose $G\subset LR(Y)$ is a semialgebra. Then for any $u\in Y
$, $G(u)\cup\{0\}$ is a linear subspace of $Y$.
\end{lemma}

\begin{proof}
Let $x,y\in G(u),$ then $x\in f(u),y\in f^{\prime}(u)$ for some $f,f^{\prime
}\in G.$ Suppose $a,b\in\mathbf{R\backslash\{0\}}$ and let $g=af+bf^{\prime
}\in G.$ Then
\[
ax+by\in af(u)+bf^{\prime}(u)=g(u)\subset G(u).
\]
\end{proof}

\begin{lemma}
\label{LinSub}Suppose $R\in LR(Y).$ Then
\[
Fix(R)=\{x\in Y:x\in R(x)\}
\]
is a linear subspace of $Y.$
\end{lemma}

\begin{proof}
Let $a,b\in\mathbf{R}\backslash\{0\},x,y\in Fix(R).$ Then $ax+by\in
aR(x)+bR(y)=R(ax+by),$ so $ax+by\in Fix(R).$
\end{proof}

\begin{lemma}
\label{FixInv}Suppose $h\in LR(Y),R\in Comm(h)$ and $h(0)\subset R(0).$ Then
$M=Fix(R)$ is $h$-weakly-invariant.
\end{lemma}

\begin{proof}
Let $x\in M=Fix(R).$ Then $A=h(x)\in hR(x)\subset Rh(x)=R(A).$ In particular,
there is some $z\in A$ such that $R(z)\cap A\neq\emptyset.$ Suppose $y\in
R(z)\cap A.$ Now we use Lemma \ref{T(0)}, as follows:
\[
z\in A=h(x)=y+h(0)\subset y+R(0)=R(z).
\]
Hence $z\in Fix(R)=M$, so $h(x)\cap M\neq\emptyset.$
\end{proof}

The main results of this paper are given below.

\begin{theorem}
[Weakly Invariant Subspaces]\label{InvSub}Suppose $Y$ is a normed space, $h\in
LR(Y)$, $G\subset Comm(h)$ is a semialgebra, $k\in Comm(h)$ is compact and
$k^{-1}(0)\neq Y.$ Suppose that for all $g\in G,$ $h(0)\subset gk(0).$ Then
there exists a nontrivial closed linear subspace $M\subset Y$ such that
\[%
\begin{tabular}
[c]{ll}%
either & (1) $M$ is $G$-invariant,\\
or & (2) $M$ is finite-dimensional and $h$-weakly-invariant.
\end{tabular}
\]
\end{theorem}

\begin{proof}
For each $u\in Y,$ $G(u)$ is $G$-invariant. Indeed, suppose $x\in G(u)$ and
$g\in G.$ Then $x\in f(u)$ for some $f\in G.$ Therefore $gf\in G$ and
$g(x)\subset gf(u)\subset G(u).$ Suppose now that $Q=G(u_{0})$ is not dense in
$Y$ for some $u_{0}\in Y\backslash\{0\}.$ Then we can assume that
$Q\neq\{0\},$ because otherwise $span\{u_{0}\}$ is $G$-invariant. Then
$M=\overline{Q}=\overline{Q\cup\{0\}}$ is the desired subspace. First,
$L=Q\cup\{0\}$ is a linear subspace of $Y$ by Lemma \ref{ComLin}. Second,
since every $f\in G$ is u.s.c., $f(Q)\subset Q$ implies $f(\overline
{Q})\subset\overline{Q}.$ Hence $M$ is $G$-invariant.

Assume now that $G(y)$ is dense in $Y$ for each $y\in Y\backslash\{0\}.$ Since
$k$ is u.s.c., $k^{-1}(0)$ is closed. Therefore we can choose a closed convex
neighborhood $U\subset Y$ of\ some $b\in Y\backslash\{0\}$ such that $0\notin
U$ and $U\subset Y\backslash k^{-1}(0).$ In addition we have
\begin{equation}
b\in\overline{G(y)},\text{ for all }y\in Y\backslash\{0\}. \tag{*}\label{inv1}%
\end{equation}
Now we apply Theorem \ref{FixNEmp} with $X=Y\backslash\{0\},$ $A=\{b\}$ as
follows. We let $H=\{g\cap U:g\in G\}.$ Then $H\subset F_{c}(Y).$ We can also
rewrite (\ref{inv1}) as
\[
\overline{H(x)}\cap A\neq\emptyset,\text{ for all }x\in Y\backslash\{0\}.
\]
By Lemma \ref{Cmpct}, there is a compact $C\subset Y$ such that $k(x)\cap
C\neq\emptyset.$ Let $r=k\cap C\in F_{c}(Y,Y)$. Since $k$ is compact, $r(U) $
is precompact. Therefore by Theorem \ref{FixNEmp}, for any neighborhood $W$ of
$b,$ $Conv(H)r$ has a fixed point in $W.$ Therefore there is $g\in
Conv(G)\subset G$ such that $R=gk$ has a fixed point $y_{0}\neq0.$ Thus
\[
M=Fix(R)\neq\{0\}.
\]
Now $R\in Comm(h)$ (Lemma \ref{CommLin}), $M$ is a linear subspace of $Y$
(Lemma \ref{LinSub}), and $M$ is $h$-invariant (Lemma \ref{FixInv}).

Suppose now that $B\subset M$ is a bounded neighborhood of $0$. Since $R$ is
compact by Theorem \ref{CompactComposition}, we have $B\subset R(B)\subset
C+R(0),$ where $C$ is compact and $R(0)$ is finite dimensional. Therefore $M$
is finite dimensional.
\end{proof}

\begin{corollary}
\label{InvSubCor}Suppose $Y$ is a normed space, $h\in LR_{0}(Y)$, $k\in
Comm(h)$ is compact and $k^{-1}(0)\neq Y.$ Then there exists a nontrivial
closed $h$-weakly-invariant subspace $M\subset Y$ .
\end{corollary}

\begin{proof}
Let
\[
G=\{P(h):P\text{ is a polynomial without constant term}\}.
\]
Then $G\subset Comm(h)$ and $G$ is a semialgebra. To check the rest of the
conditions of the theorem, observe that since $0\in k(0),$ we have
$h^{n}(0)\subset h^{n}k(0)$ for all $n\geq0.$ Therefore $h(0)\subset
h^{n}(0)\subset h^{n}k(0)$ for all $n\geq1,$ so that $h(0)\subset P(h)(0)$ for
any polynomial $P$ without constant term$.$ Thus for all $g\in G,$
$h(0)\subset gk(0).$ Therefore by the theorem there is a nontrivial closed
linear subspace $M\subset Y$ which is either $G$-invariant or $h$%
-weakly-invariant. Since $h\in G,$ $M$ is $h$-weakly-invariant.
\end{proof}

\begin{remark}
When both $h$ and $k$ are single valued, the corollary reduces to Lomonosov's
Theorem \ref{LomOrig}. The corollary is vacuous when $h$ is single-valued
while $k$ is not, because $hk\subset kh$ implies that $k(0)$ is a
finite-dimensional $h$-invariant subspace. Yet in the next section we will
obtain some applications of Theorem \ref{InvSub} to linear operators.
\end{remark}

\section{Applications to Invariant Subspaces of Linear Operators.}

In this section we generalize a well known corollary to Lomonosov's Theorem
\ref{LomOrig}.

\begin{lemma}
\label{Semi}Suppose $h,k\in LR(Y).$ Then
\[
G=\{g\in Comm(h):h(0)\subset gk(0)\}
\]
is a semialgebra.
\end{lemma}

\begin{proof}
By Lemma \ref{CommLin} $Comm(h)$ is a semialgebra, so we need only to consider
the following. (1) Suppose $h(0)\subset fk(0)$ and $h(0)\subset gk(0),$
$a,b\in\mathbf{C\backslash\{0\}}.$ Then $h(0)\subset
afk(0)+bgk(0)=(af+bg)k(0)$ because $h(0)$ is a linear subspace. Thus $G$ is
closed under nonzero linear combinations. (2) Suppose $h(0)\subset fk(0),$ $g
$ commutes with $h.$ Now since $0\in g(0),$ we have
\[
h(0)\subset hg(0)\subset gh(0)\subset gfk(0).
\]
Hence $G$ is closed under compositions.
\end{proof}

Suppose $Y$ is a space over $\mathbf{C}.$ Given $h\in LR(Y),$ its
\textit{eigenvalue} $\lambda\in\mathbf{C}$ and \textit{eigenvector} $u\neq0$
satisfy
\[
\lambda u\in h(u)
\]
\ (or $\ker(\lambda Id-h)\neq\{0\}).$ The \textit{eigenspace} of $h$
corresponding to $\lambda$ is given by
\[
E_{h\lambda}=\{u\in Y:\lambda u\in h(u)\}.
\]

\begin{lemma}
\label{EigLin}If $h\in LR(Y)$ is u.s.c. then $E_{h\lambda}$ is a closed linear subspace.
\end{lemma}

\begin{proof}
The set $E_{h\lambda}=(\lambda Id-h)^{-1}(0)$ is closed because $h$ is u.s.c..
\end{proof}

The next lemma is obvious, see \cite[Chapter 6]{Cross}.

\begin{lemma}
\label{Eigen}Let $E$ be a finite dimensional space over $\mathbf{C}.$ Then any
$h\in LR_{0}(E)$ has an eigenvector.
\end{lemma}

We say that the maps $h$,$k:Y\rightarrow Y$ \textit{strictly commute} if
\[
kh=hk.
\]

\begin{theorem}
[Invariant Subspaces]\label{InvHSub}Let $Y$ be a normed space over
$\mathbf{C,}$ $h\in LR(Y)$ is u.s.c., for any $\lambda\in\mathbf{C,}$ $\lambda
Id$ is not a selection of $h$, $k\in Comm(h)$ compact, $k^{-1}(0)\neq Y.$ If a
linear operator $f:Y\rightarrow Y$ strictly commutes with $h$ and $h(0)\subset
fk(0),$ then there is a nontrivial closed $f$-invariant subspace.
\end{theorem}

\begin{proof}
By Lemma \ref{Semi}
\[
G=\{g\in Comm(h):h(0)\subset gk(0)\}
\]
is a semialgebra. Therefore by Theorem \ref{InvSub} there is a nontrivial
closed subspace $M$ that is either (1) $G$-invariant or (2) finite-dimensional
and $h$-weakly-invariant. In case of (1) $M$ is $f$-invariant because $f\in
G.$ Consider (2). If $h^{\prime}=h\cap M:M\rightarrow M$ then $h^{\prime}$ is
a linear relation with nonempty values. Then by Lemma \ref{Eigen}, $h^{\prime
}$ has an eigenvector corresponding to some $\lambda\in\mathbf{C}$. Therefore
the eigenspace $N=E_{h\lambda}$ is nonzero, closed and not equal to the whole
$Y.$ Then for each $x\in N,$ we have
\[
\lambda f(x)=f(\lambda x)\in fh(x)=hf(x).
\]
Hence $f(x)\in N.$
\end{proof}

The following result involves only single valued operators.

\begin{corollary}
Let $Y$ be a normed space over $\mathbf{C,}$ $N$ be a finite dimensional
subspace of $Y,$ $f,h,k:Y\rightarrow Y$ bounded linear operators. Suppose for
any $\lambda\in\mathbf{C,}$ $h\neq\lambda Id$, $k$ is nonzero compact$.$
Suppose also that
\[%
\begin{tabular}
[c]{l}%
(1) $fh=hf,$\\
(2) $hk-kh\in N,$\\
(3) $h(N)\subset N.$%
\end{tabular}
\]
Then there is a nontrivial closed $f$-invariant subspace.
\end{corollary}

\begin{proof}
Apply the above theorem to the linear relation $k+N.$
\end{proof}


\begin{thebibliography}{9}
\bibitem{AAB}Y.A. Abramovich, C.D. Aliprantis, and O. Burkinshaw, The
invariant subspace problem. Some recent advances, \textit{Rend. Inst. Mat.
Univ. Trieste}, XXIX Supplemento: 3-79, 1998.

\bibitem{Cross}R. Cross, ``Multivalued Linear Operators'', Marcel Dekker, 1998.

\bibitem{En}R. Engelking, ``General Topology,'' Second Edition, Heldermann
Verlag, Berlin, 1989.

\bibitem{Gorn}L. Gorniewicz, ``Topological Fixed Point Theory of Multivalued
Mappings'', Kluwer, 1999.

\bibitem{Lom}V.I. Lomonosov, Invariant subspaces for operators commuting with
compact operators, \textit{Functional Anal. Appl.,} \textbf{7} (1973), 213-214.

\bibitem{RR}H. Radjavi and P. Rosenthal, ``Invariant Subspaces'',
Springer-Verlag, New York, 1973.

\bibitem{Sav1}P. Saveliev, A Lefschetz-type coincidence theorem, \textit{Fund.
Math.}, \textbf{162} (1999), 65-89.

\bibitem{Sim1}A. Simoni\v{c}, A construction of Lomonosov functions and
applications to the invariant subspace problem, \textit{Pac. J. Math.},
\textbf{175} (1996) 1, 257-270.

\bibitem{Sim2}A. Simoni\v{c}, An extension of Lomonosov's techniques to
non-compact operators, \textit{Proc. Amer. Math. Soc.}, \textbf{348} (1996) 3, 975-995.
\end{thebibliography}
\end{document}